\documentclass[11pt]{article}

\usepackage{a4wide}
\usepackage{amssymb}
\usepackage{amsfonts}
\usepackage{amsmath}
\input xy
\xyoption{arrow} \xyoption{matrix}

\date{}

\newtheorem{proposition}{Proposition}[section]
\newtheorem{theorem}[proposition]{Theorem}
\newtheorem{lemma}[proposition]{Lemma}

\newtheorem{corollary}[proposition]{Corollary}

\def\der{\partial }

\def\nFM0{{\nu }_{F,M_0}}
\def\nFN0{{\nu }_{F,N_0}}
\def\nGN0{{\nu }_{G,N_0}}

\def\N0{ {\bf N}_0 }

\def\t{\otimes}

\def\ra{\rightarrow}

\def\Xpm{X^{\pm }}

\def\s{\sigma}
\def\Z{\mathbb{Z}}

\def\l1{{\lambda}_1}

\def\a{\alpha}
\def\a0{ {\alpha }_0}
\def\a1{ {\alpha }_1}

\def\l{\lambda}
\def\o{\omega}

\def\nFGM0{{\nu }_{F,G,M_0}}

\def\nFN0{{\nu}_{F,N_0}}

\def\sm{{\sigma}^m}

\def\sm1{{\sigma}^{-1}}

\def\smtp1{{\sigma}^{-t+1}}

\def\o{\omega }
\def\S1{S^{-1}}

\def\Xpm1{X^{\pm 1}_1}

\def\sPM1{{\sigma }^{\pm 1}}
\def\sMP1{{\sigma }^{\mp 1 }}

\def\d{\delta}

\def\di{{\rm d.ind}}

\def\L{\Lambda}

\def\CA{{\cal A}}

\def\CD{{\cal D}}

\def\Ytm1{Y^{t-1}}
\def\Yim1{Y^{i-1}}

\def\CL{{\cal L}}
\def\CM{{\cal M}}
\def\CN{{\cal N}}

\def\Aut{{\rm Aut}}

\def\Der{{\rm Der }}
\def\ad{{\rm ad }}

\def\ker{ {\rm ker } }

\def\D{ \Delta }

\def\SL2Z{ {\rm SL}_2({\bf Z}) }

\def\CL{{\cal L}}

\def\Gp1{ G^{1 , 1 } }
\def\P11{ P^{-1 , 1 } }
\def\Pp1{ P^{1 , 1 } }

\def\nCLsr{{}^\nu\kern-2pt {\cal L}^{\sigma , \rho  }}
\def\nP{{}^\nu \kern-2pt P}
\def\nL{{}^\nu\kern-2pt L}
\def\nLL{{}^\nu\kern-2pt \Lambda}
\def\nPsr{{}^\nu\kern-2pt P^{\sigma , \rho  }}
\def\nLsr{{}^\nu\kern-2pt L^{\sigma , \rho  }}
\def\nuCL{{}^\nu\kern-2pt  {\cal L}}
\def\nCLsr{{}^\nu\kern-2pt {\cal L}^{\sigma , \rho  }}
\def\nCL1m{{}^\nu\kern-2pt {\cal L}^{-1 , 1  }}
\def\x1nu{x^\frac{1}{\nu}}
\def\xm1nu{x^{-\frac{1}{\nu}}}

\def\CN{{\cal N}}
\def\ra{\rightarrow }

\def\CB{{\cal B}}

\def\CC{ {\cal C}}

\def\nAM0{{\nu }_{{\cal A},M_0}}
\def\nAN0{{\nu }_{{\cal A},N_0}}

\def\End{ {\rm End }}
\def\Der{ {\rm Der }}

\def\ad{ {\rm ad }}

\def\derij{\partial_{{\bf i}, {\bf j}}}
\def\GL{{\rm GL}}
\def\SL{{\rm SL}}

\def\di!{\frac{\der^i}{i!}}
\def\dik!{\frac{\der^k_i}{k!}}

\def\Fp{\mathbb{F}_p}

\def\derik{\der_i^{[k]}}
\def\derij{\der_i^{[j]}}

\def\deril{\der_i^{[l]}}
\def\derjl{\der_j^{[l]}}
\def\derba{\der^{[\alpha ]}}
\def\derbb{\der^{[\beta ]}}

\def\Nn{\mathbb{N}^n}

\def\id{{\rm id}}

\def\Zp{\mathbb{Z}_p}

\def\N{\mathbb{N}}

\def\0{\overline{0}}
\def\1{\overline{1}}

\def\Ln1{\L_{n,\overline{1}}}

\def\a1{a_{\overline{1}}}

\def\St{{\rm St}}
\def\S{\Sigma}

\def\vn1{\overrightarrow{n-1}}

\def\deripk{\der_i^{[p^k]}}

\def\derijkpk{\der_i^{[j_kp^k]}}
\def\deripkjk{\der_i^{[p^k]j_k}}
\def\CDPn{\CD (P_n)}

\begin{document}

\author{V. V. \  Bavula 
}

\title{The group of order preserving automorphisms of the ring of differential operators on Laurent polynomial algebra in prime characteristic }

\maketitle
\begin{abstract}
Let $K$ be a  field of characteristic $p>0$. It is proved that the
group $\Aut_{ord}(\CD (L_n))$ of order preserving automorphisms of
the ring $\CD (L_n)$ of differential operators on a Laurent
polynomial algebra $L_n:= K[x_1^{\pm 1}, \ldots , x_n^{\pm 1}]$ is
isomorphic to a skew direct product of groups $\Zp^n \rtimes
\Aut_K(L_n)$ where $\Zp$ is the ring of $p$-adic integers.
Moreover, the group $\Aut_{ord}(\CD (L_n))$ is found explicitly.
Similarly, $\Aut_{ord}(\CDPn )\simeq \Aut_K(P_n)$ where $P_n:
=K[x_1, \ldots , x_n]$ is a polynomial algebra.



\end{abstract}


\section{Introduction}
Throughout, ring means an associative ring with $1$; $\N :=\{0, 1,
\ldots \}$ is the set of natural numbers; $p$ is a prime number;
$\Fp := \Z / \Z p$ is the field that contains $p$ elements; $\Zp$
is the ring of $p$-adic integers;
 $K$ is an arbitrary field of characteristic $p>0$ (if it is not
stated otherwise); $P_n:= K[x_1, \ldots , x_n]$ is a polynomial
algebra;
 $\CDPn =\oplus_{\alpha \in \N^n}P_n\der^{[\alpha]}$ is the ring of
differential operators on $P_n$ where
$\der^{[\alpha]}:=\prod_{i=1}^n
\frac{\der_i^{\alpha_i}}{\alpha_i!}$; $\D_n:=\oplus_{\alpha \in
\N^n}K\der^{[\alpha]}$ is the algebra of scalar differential
operators on $P_n$; $L_n:=K[x_1^{\pm 1}, \ldots , x_n^{\pm 1}]$ is
a Laurent polynomial algebra and $\CD (L_n) =\oplus_{\alpha \in
\N^n}L_n\der^{[\alpha]}$ is the ring of differential operators on
the algebra  $L_n$; $\{ \CD (L_n)_i:= \oplus_{|\alpha |\leq i}
L_n\der^{[\alpha]}\}_{i\geq 0}$ is the {\em order filtration} on
$\CD (L_n)$;
$$ \Aut_{ord} (\CD (L_n)):=\{ \s \in \Aut_K(\CD (L_n))\, | \; \s
(\CD (L_n)_i)=\CD (L_n)_i, \; i\geq 0\}$$ is the group of {\em
order preserving} automorphisms of the algebra $\CD (L_n)$;
similarly the order filtration on $\CDPn$ and the group
$\Aut_{ord}(\CDPn )$ are defined.

In arbitrary characteristic, it is a difficult problem to find
generators for the groups $\Aut_K(P_n)$ and $\Aut_K(\CDPn )$. The
results are known for $\Aut_K(P_1)$ (easy), $\Aut_K(P_2)$ (Jung
\cite{Jung} and van der Kulk \cite{vderKulk}), and $\Aut_K(\CD
(P_1))$ when ${\rm char}(K)=0$ (Dixmier \cite{Dix-1968}). Little
is known about the groups of automorphisms in the remaining cases.

 In
characteristic zero, there is a strong connection between the
groups $\Aut_K(P_n)$ and $\Aut_K(\CDPn )$ as, for example, the
(essential) equivalence of the Jacobian Conjecture for $P_n$ and
the Dixmier Problem/Conjecture for $\CDPn$ shows (see \cite{BCW},
\cite{Tsuchi05}, \cite{Bel-Kon05JCDP}, see also  \cite{JC-DP}).
Moreover,  conjecture like the two mentioned conjectures makes
sense only for the algebras $P_m\t \CDPn$ as was proved in
\cite{Bav-cpinv} (the two conjectures can be reformulated in terms
of locally nilpotent derivations that satisfy certain conditions,
and the algebras $P_m\t \CDPn$  are the only associative algebras
that have such derivations). This general conjecture is true iff
either the JC  or the DC is true, see \cite{Bav-cpinv}.

In prime characteristic,  relations between the two groups,
$\Aut_K(P_n)$ and $\Aut_K(\CDPn )$, are  even tighter as the
following  result shows.

\begin{theorem}\label{28Jun7}
\cite{Bav-Frobext} {\rm (Rigidity of the group $\Aut_K(\CDPn )$)}
Let $K$ be a field of characteristic $p>0$, and $\s , \tau \in
\Aut_K(\CDPn )$. Then $\s = \tau $ iff $\s (x_1) = \tau (x_1) ,
\ldots , \s (x_n) = \tau (x_n)$.
\end{theorem}

{\it Remark}.  Theorem \ref{28Jun7} does not hold in
characteristic zero, and, in general,  in prime characteristic it
does not hold for localizations of the polynomial algebra $P_n$
(Theorem \ref{10May8}). Note that the $K$-algebra $\CDPn$ is not
finitely generated.

$\noindent $

As a direct consequence of Theorem \ref{28Jun7} there is the
corollary (see Section \ref{TGAD} for details).

\begin{corollary}\label{aa28Jun8}
$\Aut_{ord}(\CDPn ) \simeq \Aut_K(P_n)$.
\end{corollary}

Situation is completely different for the Laurent polynomial
algebra $L_n$.

\begin{theorem}\label{10May8}
$\Aut_{ord}(\CD (L_n))\simeq \Zp^n \rtimes \Aut_K(L_n)$. Moreover,
each automorphism $\s \in\Aut_{ord}(\CD (L_n))$ is a unique
product of two automorphisms $\s = \s_s\tau$ where $\tau \in
\Aut_K(L_n)$, $s=(s_i)\in \Zp^n$, $\s _s (x_i) = x_i$ for all $i$,
and 
\begin{equation}\label{ssad}
\s_s(\der^{[\alpha]})= \prod_{i=1}^n
\frac{(\der_i+s_ix_i^{-1})^{\alpha_i}}{\alpha_i!}, \;\; \alpha
=(\alpha_i)\in \N^n.
\end{equation}
\end{theorem}

The meaning of the RHS of (\ref{ssad}) is explained in Section
\ref{TGAD}. Theorem \ref{10May8} gives all the elements of the
group $\Aut_{ord}(\CD (L_n))$ explicitly since the group
$\Aut_K(L_n)$ is known, it is isomorphic to the semidirect product
$\GL_n(\Z )\rtimes K^{*n}$.

 The main idea of the proof of Theorem
\ref{10May8} is to show that the group $\Aut_{ord} (\CD (L_n))$ is
equal to the group
$$\St (L_n) :=\{ \s \in \Aut_K(\CD (L_n))\, | \; \s (L_n) =
L_n\}$$ which is a semidirect product $\St (L_n)={\rm st} (L_n)
\rtimes \Aut_K(L_n)$ where
$${\rm st}  (L_n) :=\{ \s \in \St (L_n)\, | \; \s |_{L_n} =
\id \}$$ is the normal subgroup of $\St (L_n)$. The group ${\rm
st}(L_n)$ consists of all the automorphisms $\s_s$, $s\in \Zp^n$
(see (\ref{ssad})), and the map $\Zp^n\ra {\rm st}(L_n)$,
$s\mapsto \s_s$, is a group isomorphism ($\s_{s+t} = \s_s\s_t$),
Theorem \ref{31Mar8}.


\section{The group $\Aut_{ord}(\CD (L_n))$} \label{TGAD}

Let $R$ be a commutative $K$-algebra and $\CD (R)=\cup_{i\geq
0}\CD (R)_i$ be the ring of $K$-linear differential operators on
the algebra $R$ where $\{ \CD (R)_i\}_{i\geq 0}$ is the {\em order
filtration} on $\CD (R)$. In particular, $\CD (R)_0= R$. Let
$\Aut_{ord} (\CD (R))$ be the subgroup of $\Aut_K(\CD (R))$ of
order preserving $K$-automorphisms, i.e.
$$\Aut_{ord}(\CD (R)):= \{ \s \in \Aut_K(\CD (R))\, | \; \s (\CD
(R)_i)= \CD (R)_i,\;  i\geq 0\}.$$ Each automorphism $\s \in
\Aut_K(R)$ can be naturally extended (by change of variables) to a
$K$-automorphism, say $\s$, of the ring $\CD (R)$ of differential
operators on the algebra $R$ by the rule: 
\begin{equation}\label{hs=sas}
\s(a) := \s a \s^{-1}, \;\; a\in \CD (R) .
\end{equation}
Then the group $\Aut_K(R)$ can be seen as a subgroup of $\Aut_K
(\CD (R) )$ via (\ref{hs=sas}). Then it follows from  a definition
of the ring $\CD (R)$ and the equality $[r, \s a \s^{-1} ] = \s
[\s^{-1} (r), a]\s^{-1}$ where $r\in R$, $a\in \CD (R)$ that
\begin{equation}\label{1hs=sas}
\Aut_K(R)\subseteq \Aut_{ord}(\CD (R)).
\end{equation}

The stabilizer of the algebra $R$,
$$ \St (R):= \{ \s \in \Aut_K(\CD (R) )\, | \,  \s (R) = R\}$$
is a subgroup of $\Aut_K(\CD (R) )$. It contains the normal
subgroup
$${\rm st} (R) := \{ \s \in \St (R) \, | \, \s|_{R}=\id\}$$
which is the kernel of the restriction epimorphism
$$\St (R) \ra
\Aut_K(R), \;\; \s\mapsto \s|_{R}.$$
 By (\ref{1hs=sas}),  the
 stabilizer of $R$, 
\begin{equation}\label{StRsd}
\St (R) = {\rm st}(R) \rtimes \Aut_K(R)
\end{equation}
is the semi-direct product of its subgroups. It is obvious that
\begin{equation}\label{AiSt}
\Aut_{ord}(\CD (R))\subseteq \St (R).
\end{equation}
The ring $R$ is a left $\CD (R)$-module. The action of an element
$\d\in \CD (R)$ on an element $r\in R$ is denoted either by $\d
(a)$ or $\d *r$ (in order to avoid multiple brackets).

$\noindent $

{\bf Proof of Corollary \ref{aa28Jun8}}. By Theorem \ref{28Jun7},
${\rm st}(P_n) = \{ \id \}$. Now, the result follows from
(\ref{1hs=sas}), (\ref{StRsd}) and (\ref{AiSt}):

$$  \Aut_K(P_n) \subseteq \Aut_{ord}(\CDPn )\subseteq \St (P_n) = {\rm st}(P_n) \rtimes \Aut_K(P_n)  =\Aut_K(P_n).\;\; \Box$$

$\noindent $

{\bf The rings  of  differential operators $\CD (P_n)$   and $\CD
(L_n)$}. The ring $\CD (P_n)$ of differential operators on a
polynomial algebra $P_n:= K[x_1, \ldots , x_n]$ is a $K$-algebra
generated by the elements $x_1, \ldots , x_n$ and {\em commuting}
higher derivations $\derik :=\frac{\der_i^k}{k!}$, $i=1, \ldots ,
n$; $k\geq 1$,  that satisfy the following defining relations:
\begin{equation}\label{DPndef}
[x_i,x_j]=0, \;\; [\derik , \derjl ]=0,\;\;\; \derik \deril
={k+l\choose k}\der_i^{[k+l]}, \;\;\; [\derik ,
x_j]=\d_{ij}\der_i^{[k-1]},
\end{equation}
 for all
$i,j=1, \ldots , n$;  $k,l\geq 1$, where $\d_{ij}$ is the
Kronecker delta,  $\der_i^{[0]}:=1$, $\der_i^{[-1]}:=0$, and
$\der_i^{[1]}=\der_i=\frac{\der}{\der x_i}\in \Der_K(P_n)$,
$i=1,\ldots , n$. The action of the higher derivation $\derik$ on
 the polynomial algebra
 $$P_n=K\t_\Z \Z [x_1, \ldots , x_n]\simeq K\t_\Z \Zp [x_1, \ldots , x_n]$$
should
 be understood as the action of the element $1\t_\Z
 \frac{\der_i^k}{k!}$.

The algebra $\CD (P_n)$ is a {\em simple} algebra. Note  that the
algebra $\CD (P_n)$ is  not finitely generated and not (left or
right) Noetherian, it  does not satisfy finitely many defining
relations.
$$\CD (P_n)= \bigoplus_{\alpha , \beta \in \Nn }Kx^\alpha \derbb
\subset \CD (L_n) = \bigoplus_{\alpha\in \Z^n , \beta \in \Nn }K
 x^\alpha \derbb
 .$$
For each $i=1, \ldots , n$ and $j\in \N$ written $p$-adically as
$j= \sum_k j_kp^k$, $0\leq j_k<p$, 
\begin{equation}\label{Fdij}
\derij=\prod_k\derijkpk = \prod_k \frac{\deripkjk}{j_k!}, \;\;
\derijkpk =\frac{\deripkjk}{j_k!},
\end{equation}
where $\deripkjk := (\deripk )^{j_k}$. For $\alpha   \in \Nn$ and
$\beta \in \Z^n$, 
\begin{equation}\label{daaxb}
\derba *x^\beta  = {\beta \choose \alpha } x^{\beta - \alpha},
\;\; {\beta \choose \alpha }:= \prod_i {\beta_i \choose \alpha_i}.
\end{equation}

For $\alpha , \beta \in \Nn$, 
\begin{equation}\label{dadb}
\derba \derbb = {\alpha +\beta \choose \beta}\der^{[\alpha +\beta
]}.
\end{equation}

{\bf The binomial differential operators}.
 For each
natural number $i$, the binomial polynomial
$$ {t\choose i} :=\frac{t(t-1)\cdots (t-i+1)}{i!}, \;\; {t\choose
0}:=1,$$ can be seen as a function from $\Zp$ to $\Zp$. For each
$p$-adic integer  $s\in \Z$ and a natural number  $i\in \N$, we
have the differential operator on the $\Zp$-algebra $\Zp [x^{\pm
1}]:= \Zp [x,x^{-1}]$ of Laurent polynomials with coefficients
from $\Zp$: 
\begin{equation}\label{pds1x}
\frac{(\der +sx^{-1})^i}{i!}=x^{-i} \frac{(x\der +s) (x\der +s-1)
\cdots ( x\der +s-i+1)}{i!}=x^{-i}{x\der +s\choose i}\in \CD (\Zp
[x^{\pm 1}]),
\end{equation}
where $\der := \frac{d}{dx}\in \Der_{\Zp} (\Zp [x^{\pm 1}])$. Let
$\CL_n:= \Zp [ x_1^{\pm 1}, \ldots , x_n^{\pm 1}]$ be the Laurent
polynomial ring in $n$ variables with coefficients from $\Zp$ and
let  $\CD (\CL_n)$ be the ring of $\Zp$-linear differential
operators on the ring $\CL_n$. Then, for any elements $\alpha
=(\alpha_i)\in \N^n$ and $s=(s_i)\in \Zp^n$, there is  the
differential operator  on the algebra $\CL_n$: 
\begin{equation}\label{p1ds1x}
 b_s^{[\alpha ]}:=\prod_{i=1}^n x_i^{-\alpha_i}{x_i\der_i +s_i\choose \alpha_i}=
 \prod_{i=1}^n\frac{(\der_i +s_ix_i^{-1})^{\alpha_i}}{\alpha_i!}  \in \CD
 (\CL_n).
\end{equation}

The inclusions of abelian monoids $\N^n\subset \Z^n \subset \Zp^n$
yield the inclusions of their monoid rings
$$ \CN_n:= \bigoplus_{\alpha \in \N^n} \Zp x^\alpha \subset
\CL_n\subset \CM_n:= \bigoplus_{\alpha \in \Zp^n} \Zp x^\alpha ,
$$ and the inclusions of rings
$$\CA_n:=\bigoplus_{\beta
\in \N^n}\CN_n \der^{[\beta]}  \subset \CB_n:= \bigoplus_{\beta
\in \N^n}\CL_n \der^{[\beta]}  \subset \CC_n:= \bigoplus_{\beta
\in \N^n}\CM_n \der^{[\beta]}\subset \End_{\Z_p} (\CM_n)$$ where
$\der^{[\beta]}= \prod_{i=1}^n \frac{\der_i^{\beta_i}}{\beta_i!}$.
For any elements $\alpha \in \Zp^n$ and $\beta \in \N^n$,
$$ \der^{[\beta]}*x^\alpha = {\beta \choose \alpha} x^{\alpha -
\beta}\;\; {\rm where}\;\; {\beta \choose \alpha}:= \prod_{i=1}^n
{\beta_i \choose \alpha_i}.$$ The prime number $p$ belongs to the
centre of each of the rings $\CA_n$ , $\CB_n$ and $\CC_n$. Taking
factor-rings modulo the ideals generated by the element $p$ in
each of the rings we obtain the inclusions of $\Fp$-algebras:
$$ \bigoplus_{\alpha , \beta \in \N^n}\Fp x^\alpha \der^{[\beta]}\subset  \bigoplus_{\alpha\in \Z^n,\beta \in \N^n }\Fp x^\alpha
\der^{[\beta]}\subset \bigoplus_{\alpha\in \Zp^n, \beta \in
\N^n}\Fp x^\alpha \der^{[\beta]}\subset \End_{\Fp}
(\bigoplus_{\alpha\in \Z^n}\Fp x^\alpha ).$$ Then, applying
$K\t_\Z -$ we obtain the inclusions of rings
$$ \CDPn \subset \CD (L_n)\subset \CD_n:=\bigoplus_{\alpha \in
\Zp^n, \beta \in \N^n} Kx^\alpha \der^{[\beta]}.$$

For each element $s\in \Zp^n$, the inner automorphism
$\o_{x^{-s}}: a\mapsto x^{-s}ax^s$ of the ring $\CC_n$ acts
trivially  on the ring $\CM_n$, and,  for each element $\beta \in
\N^n$, 
\begin{equation}\label{wxms}
\o_{x^{-s}}(\der^{[\beta]})= b_s^{[\beta]}.
\end{equation}
Indeed,
\begin{eqnarray*}
 \o_{x^{-s}}(\der^{[\beta]})&= & \prod_{i=1}^n x_i^{-s_i}\der_i^{[\beta_i]}x_i^{s_i}=
 \prod_{i=1}^n x_i^{-s_i}x_i^{-\beta_i}\frac{x_i^{\beta_i}\der_i^{\beta_i}}{\beta_i!}x_i^{s_i}=
 \prod_{i=1}^n x_i^{-\beta_i}x_i^{-s_i}{x_i\der_i\choose \beta_i} x_i^{s_i}\\
 & =& \prod_{i=1}^n x_i^{-\beta_i}{x_i\der_i+s_i\choose \beta_i}x_i^{-s_i}
 x_i^{s_i}=b_s^{[\beta]}.
\end{eqnarray*}
The inner automorphism $\o_{x^{-s}}: b\mapsto x^{-s}bx^s$ of the
ring $\CD_n$ is the restriction modulo $p$ of the inner
automorphism $\o_{x^{-s}}: a\mapsto x^{-s}ax^s$ of the ring
$\CC_n$. It is obvious that
$$ \o_{x^{-s}}(\CD (L_n))= \CD(L_n). $$
Let $\s_s$ be the restriction of the inner automorphism
$\o_{x^{-s}}$ of the ring $\CD_n$ to  the subring $\CD (L_n)$.
Then the automorphism $\s_s$ is as in (\ref{ssad}) and $\s_s\in
{\rm st}(L_n)$.

It is obvious that the map 
\begin{equation}\label{Zpst}
\Zp^n\ra {\rm st}(L_n), \;\; s\mapsto \s_s,
\end{equation}
is a group monomorphism. In fact, this map is an isomorphism,
Theorem \ref{31Mar8}. Before giving the proof of Theorem
\ref{31Mar8}, we need to prove some more results.

\begin{lemma}\label{a28Mar8}
The kernel of the $\Fp$-linear map $\der_i^{p-1}+F$ acting on the
field $L_n$ is $\Fp x_i^{-1}$ where $F: a\mapsto a^p$ is the
Frobenius on $L_n$.
\end{lemma}

{\it Proof}. First, let us consider the case when $i=n=1$. We drop
the subscript 1 in this case. Let $s= \l x^j+\cdots \in K[x^{\pm
1}]$ be a nonzero element of the kernel of the $\Fp$-linear map
$\der^{p-1}+F$ where $\l x^j$ is its least term, $0\neq \l \in K$.
Then
$$ 0=(\der^{p-1}+F)(\l  x^j) = \l j (j-1) \cdots (j-p+2)
x^{j-p+1} + \l^p x^{pj},$$ and so $j-p+1=pj$ and $ \l j (j-1)
\cdots (j-p+2) +\l^p=0$. The first equality yields $j=-1$, then
the second can be written as $\l^p-\l =0$ since $(p-1)!\equiv -1
\mod p$. Therefore, $\l \in \Fp$ and the set $\Fp x^{-1}$ belongs
to the kernel of the map $\der^{p-1}+F$. The least term of the
element $s-\l x^{-1}$ of the kernel of the map  $\der^{p-1}+F$ has
to be zero by the above argument, and so $s -\l x^{-1}=0$. This
proves that $\ker (\der^{p-1}+F)= \Fp x^{-1}$. Now, the general
case follows from this special one since $L_n\subset K(x_1, \ldots
, x_{i-1}, x_{i+1}, \ldots ,x_n)[ x_i^{\pm 1}]$. $\Box $

$\noindent $

Let $R$ be a ring. The {\it first Weyl algebra} $A_1(R)$ over $R$
is a ring generated over $R$ by two elements $x$ and $\der$ that
satisfy the defining relation $\der x - x\der =1$.

\begin{theorem}\label{17Jul7}
\cite{Bav-A1rescen}. Let $K$ be a reduced commutative
$\Fp$-algebra and $A_1(K)$ be the first Weyl algebra over $K$.
Then $$ (\der +f)^p= \der^p+\frac{d^{p-1}f}{dx^{p-1}}+f^p$$ for
all $f\in K[x]$. In more detail, $(\der + f)^p= \der^p-\l_{p-1}
+f^p$ where $f= \sum_{i=0}^{p-1} \l_i x^i\in K[x]=
\oplus_{i=0}^{p-1} K[x^p]x^i$, $\l_i\in K[x^p]$.
\end{theorem}

\begin{corollary}\label{a4Apr8}
For each element $f\in L_n$ and $i=1, \dots, n$, there is the
equality
$$ (\der_i+f)^p = \frac{\der^{p-1} f}{\der x_i^{p-1}}+f^p, $$
in the ring $\CD (L_n)$.
\end{corollary}

{\it Proof}. Let $L$ be denote the  subalgebra $K[x_1^{\pm 1},
\ldots , x_{i-1}^{\pm 1}, x_i^{\pm p},  x_{i+1}^{\pm 1}, \ldots ,
x_n^{\pm 1}]$ of the algebra $L_n$. Then $L_n
=\oplus_{j=0}^{p-1}Lx_i^j$. Consider the $L$-algebra homomorphism:
$$ A_1(L):= L\langle x, \der \rangle \ra \CD (L_n),\;\;  x\mapsto
x_i, \;\; \der\mapsto \der_i, \;\; l\mapsto l, $$ where $l\in L$.
Now, the result follows from Theorem \ref{17Jul7} since
$\der_i^p=0$. $\Box $

\begin{theorem}\label{31Mar8}
The map (\ref{Zpst}) is an isomorphism.
\end{theorem}

{\it Proof}. It suffices to show that, for  each automorphism $\s
\in {\rm st}(L_n)$,
$$ \cdots \s_{p^ks_k}\cdots \s_{ps_1}\s_{s_0}\s = 1$$
for some elements $s_k= (s_{k1}, \ldots , s_{kn})\in \{ 0,1,
\ldots , p-1\}^n$. Note that, for all $k\geq 0$ and $i=1, \ldots ,
n$, 
\begin{equation}\label{spksk}
\s_{p^ks_k}(\der_i^{[p^s]}) = \der_i^{[p^s]}, \;\; s<k; \;\;
\s_{p^ks_k}(\der_i^{[p^k]}) = \der_i^{[p^sk]}+s_{ki}x_i^{-p^k}.
\end{equation}
For all indices $i$ and $j$,
$$ [ \s (\der_i)-\der_i, x_j]=\s ([\der_i, x_j])-[\der_i, x_j]=
\d_{ij}-\d_{ij}=0,$$ hence $a_i:= \s (\der_i) -\der_i \in C(\CD
(L_n) , x_1, \ldots , x_n):= \bigcap_{i=1}^n\ker (\ad (x_i))=
L_n$, the centralizer of the elements $x_1, \ldots , x_n$ in the
ring $\CD (L_n)$ where $\ad (x_i)$ is the inner derivation of the
ring $\CD (L_n)$  determined by the element $x_i$. Using Corollary
\ref{a4Apr8}, we see that
$$ 0=\s(\der_i^p)=(\der_i+a_i)^p=(\der_i^{p-1}+F)(a_i),$$
and so $a_i= -s_{0i}x_i^{-1}$ for some elements $s_{0i}\in \Fp$,
by Lemma \ref{a28Mar8}. Abusing the notation we assume in this
proof that $\Fp = \{ 0, 1, \ldots , p-1\}$.
 Let $s_0:= (s_{01}, \ldots s_{0n})$. Then
$\s_{s_0}\s (\der_i) = \der_i$ for all $i$. Suppose that $k>0$ and
we have found vectors $s_0, \ldots , s_{k-1} \in \Fp^n$ such that
$$ \s_{p^{k-1}s_{k-1}}\cdots\s_{ps_1} \s_{s_0}\s (\der_i^{[p^s]})=
\der_i^{[p^s]}, \;\; 1\leq i \leq n , \;\; 0\leq s \leq k-1.$$ In
order to finish the proof of the claim by induction on $k$, we
have to find a vector $s_k\in \Fp^n$ such that the above
equalities hold for $k$ rather than $k-1$. Let $\tau :=\s_{p^{k-1}
s_{k-1}}\cdots \s_{ps_1}\s_{s_0}\s$. For all $i$ and $j$,
$$ [ \tau (\deripk ) - \deripk, x_j^{p^k}]= \tau ([\deripk ,
x_j^{p^k}]) - [ \deripk , x_j^{p^k}]= \d_{ij} - \d_{ij} =0,$$
hence $b_i:= \tau (\deripk ) - \deripk \in C(\CD (L_n) ,
x_1^{p^k}, \ldots x_n^{p^k})= L_n\t \D_{n, <p^k}$ where $\D_{n,
<p^k}=\oplus_{{\rm all}\; \alpha_i<p^k}K\der^{[\alpha ]}$. For all
$i\neq j$,
$$ 0=\tau ([\deripk , x_j]) = [ \deripk+b_i, x_j]=[b_i, x_j],$$
and so $b_i\in \sum_{s<p^k} L_n\der_i^{[s]}$. Now,
$$ \der_i^{[p^k-1]}= \tau ( \der_i^{[p^k-1]})=\tau ([ \deripk ,
x_i]) = [\deripk +b_i, x_i] = \der_i^{[p^k-1]} + [b_i , x_i],$$
hence $b_i\in L_n$. For all $i$ and $j$, and $s<k$,
$$ 0=\tau([ \deripk , \der_j^{[p^s]}]) = [ \deripk +b_i, \der_j^{[p^s]}]= [b_i, \der_j^{[p^s]}], $$
and so $b_i\in K[x_1^{\pm p^k}, \ldots , x_n^{\pm p^k}]$. Finally,
$$ 0=\tau (\der_i^{[p^k]p} ) = (\deripk +b_i)^p = (
\der_i^{[p^k](p-1)}+F)(b_i),
$$ hence $b_i= -s_{ki}x_i^{p^k}$ for some elements $s_{ki}\in \Fp$, by Lemma \ref{a28Mar8}. Let
$s_k= (s_{k1}, \ldots , s_{kn})$. Then
$$ \s_{p^ks_k} \cdots\s_{ps_1} \s_{s_0}\s (\der_i^{[p^s]}) =
\der_i^{[p^s]}, \;\; 1\leq i\leq n , \;\; 0\leq s\leq k.$$ This
finishes the proof of the theorem.  $\Box $

\begin{corollary}\label{a31Mar8}
 $\Aut_{ord}(\CD
(L_n) ) = \St (L_n)$.
\end{corollary}

{\it Proof}. The  inclusion $\Aut_{ord}(\CD (L_n) )\subseteq \St
(L_n)$ is obvious, see (\ref{AiSt}). The inverse inclusion follows
from the facts that $\St (L_n) = {\rm st}(L_n) \rtimes
\Aut_K(L_n)$, ${\rm st}(L_n) \subseteq \Aut_{ord} (\CD (L_n) )$
(by Theorem \ref{31Mar8}) and $\Aut_K(L_n) \subseteq
\Aut_{ord}(\CD (L_n) )$, see (\ref{1hs=sas}). $\Box $

$\noindent $

{\bf Proof of Theorem \ref{10May8}}. This follows from Corollary
\ref{a31Mar8}, (\ref{StRsd}) and Theorem \ref{31Mar8}:
 $$ \Aut_{ord} (\CD (L_n)) = \St (L_n) = {\rm st}(L_n) \rtimes \Aut_K(L_n) \simeq \Zp^n \rtimes \Aut_K(L_n).   \;\;\; \Box $$



Department of Pure Mathematics

University of Sheffield

Hicks Building

Sheffield S3 7RH

UK

email: v.bavula@sheffield.ac.uk







%

\end{document}